\documentclass[a4 paper,10pt]{article}
\usepackage{amsthm}
\usepackage{amsfonts}
\usepackage{amsmath}
\usepackage{amssymb}
\newtheorem{defi}{Definition}
\newtheorem{teo}[defi]{Theorem}
\newtheorem*{rem}{Remark}
\newtheorem{prop}[defi]{Proposition}

\newtheorem{lemma}[defi]{Lemma}
\newtheorem{ex}{Example}

\newcommand{\eps}{\varepsilon}
\newcommand{\R}{I\!\!R}
\newcommand{\ox}{\overline{X}}
\newcommand{\oa}{\overline{\alpha}}
\newcommand{\lt}{\tilde{l}}
\begin{document}
\title{Almost sure stabilizability of controlled degenerate diffusions
\thanks{This research was partially supported by M.I.U.R.,
project ``Viscosity, metric, and control theoretic methods for nonlinear
partial differential equations'', and by GNAMPA-INDAM, project ``Partial
differential equations and control theory''.}} \author{Martino
Bardi \quad and  \quad Annalisa Cesaroni\\
Dipartimento di Matematica P.  e A.\\ Universit\`a di Padova\\
via Belzoni 7, 35131 Padova, Italy\\ bardi@math.unipd.it \quad and \quad
acesar@math.unipd.it\\ \\
}

\date{}
 \maketitle
\begin{abstract}
We develop a direct Lyapunov
method for the almost sure open-loop stabilizability and asymptotic
stabilizability of controlled degenerate diffusion processes.
The infinitesimal decrease condition for a  Lyapunov function is a new
form of Hamilton-Ja\-co\-bi-Bell\-man partial differential inequality
of $2nd$ order.  We give local and global versions of the First and
Second Lyapunov Theorems assuming the existence of a lower
semicontinuous Lyapunov function satisfying such inequality in the
viscosity sense.  An explicit formula for a stabilizing feedback is
provided for affine systems with smooth Lyapunov function.  Several
examples illustrate the theory.


\end{abstract}

\section{Introduction}
For controlled diffusion processes in $\R^{N}$
\[(CSDE)\left\{\begin{array}{l}
dX_t=f(X_t,\alpha_t)dt
+\sigma(X_t,\alpha_t)dB_t,\;\;\alpha_t\in A,\;\; t>0,\\ X_0=x,
\end{array}
\right. \]
there are various possible notions of Lyapunov stability of an
equilibrium, say the origin. The stability in probability has been
studied for a long time, we recall here the contributions of  Kushner
\cite{kus, kus2},  Has'minskii \cite{has}, and the recent book of Mao
\cite{maolib} for uncontrolled systems, and the work of Florchinger
\cite{flo, flo2, flo3} and Deng, Krsti\'c, and Williams \cite {deng} on
feedback stabilization for $(CSDE)$, see also the references therein.
The almost sure exponential stability was introduced and studied by
Kozin \cite{koz}, see also \cite{has}, and it implies that for each
fixed sample in a set of probability 1 the (uncontrolled) system is
exponentially stable in the usual sense.
In this paper we consider a property that we call almost
sure stability , or uniform stability with probability 1.  For an
uncontrolled system it says that for any $\eta>0$ there exists
$\delta>0$ such that, for any 
$x$ with $|x|\leq \delta$, the process satisfies $|X_t|\leq \eta$
for all $t\geq 0$ almost surely.
Equivalently, for some increasing, continuous function $\gamma$ null at 0,
and for small $|x|$,
\begin{equation}
\label{cappa0}
|{X}_t|\leq
\gamma(|x|)\quad \forall t\geq 0 \,\text{ a.s.}
\end {equation}
This property describes a behaviour very similar to a stable
deterministic system, it is stronger than stability in probability and
pathwise stability, and in fact it is never verified by a nondegenerate
process.  More precisely, we study the {\em almost sure open loop
stabilizability} of $(CSDE)$, namely, that for each $x$ as above there
exists an admissible control function
whose trajectory $\overline{X}_{.}$ verifies a.s. $|\overline{X}_t|\leq
\eta$ (and $|\ox_t|\leq\gamma(|x|)$) for all $t$. If, in addition,
$\lim_{t\to +\infty}\overline{X}_t=0$ a.s., we say the system is {\em a.s. (open
loop) asymptotically stabilizable}. For deterministic systems
($\sigma\equiv 0$) the last
property reduces to the well-known {\em asymptotic controllability}.

We follow the Lyapunov Direct Method and find that the
{\em infinitesimal decrease condition} to be satisfied by a Lyapunov function
$V$ for our problem is
\begin{equation}
\label{cseq1.0}
\max_{\alpha\in A,\, \sigma(x,\alpha)^T DV(x) = 0}
\left\{-DV(x)\cdot
f(x,\alpha)-trace\left[a(x,\alpha) D^2V(x)\right]
\right\} \geq l(x) , 
\end{equation}
with $l\geq 0$ for mere Lyapunov stability and $l>0$ for $x\ne 0$ for
asymptotic stability, where $a:=\sigma\sigma^T/2$.
This is not a standard Hamilton-Jacobi-Bellman
inequality, because the constraint on the control $\alpha$ depends on
$V$.  In fact it should be rather viewed as a system of PDEs and
inequalities, that in the special case of uncontrolled diffusion,
i.e., $\sigma=\sigma(x)$, reads
\begin{equation}
\label{cseq1.1}
\left\{\begin{array}{l}
 \max_{\alpha\in A}\left\{-DV(x)\cdot
f(x,\alpha)\right\}-trace\left[a(x) D^2V(x)\right] \geq l(x)
 \\ \sigma_i(x)\cdot DV(x)=0 \quad \forall i ,
 \end{array}\right.
\end{equation}
where $\sigma_i$ denotes the $i$-th column of the matrix
$\sigma$.
Note that the first condition of the system is the
Hamilton-Jacobi-Bellman inequality that defines the Lyapunov functions
for the stability in probability.  The other conditions mean that
there is diffusion only in the directions tangential to the level sets
of $V$.
To motivate the infinitesimal decrease condition \eqref{cseq1.1} we
take $V$ of class $C^2$ and suppose that for a fixed control
function 
$dV(X_t)/dt\leq l(X_t)$. Then Ito's formula gives
\begin{multline*}
        \left[DV(X_t)\cdot f(X_t,\alpha_t)+
trace\left(a(X_t) D^2V(X_t)\right)\right]dt+
\sigma^T(X_t) DV(X_t)
dB_t\\
\leq l(X_t) .
\end{multline*}
Now the properties of the Brownian motion
lead to the conditions
\begin{eqnarray*}
DV(X_t)\cdot f(X_t,\alpha_t)+trace\left(a(X_t) D^2V(X_t) \right)
& \leq & l(X_t) ,\\
\sigma^T(X_t) DV(X_t) & = & 0 ,
\end{eqnarray*}
and the existence of a control $\alpha_t$ verifying this is clearly
related to \eqref{cseq1.1}.

We define a Lyapunov function for the a.s. stability as a {\em lower
semicontinuous} proper function $V$, continuous at 0 and satisfying
\eqref{cseq1.0}  in the {\em viscosity sense}, and we call it strict Lyapunov
function if $l>0$ off 0, see the Definitions
\ref{cliap} and  \ref{scliap} below. Our main results are the natural
extensions to the controlled diffusions of the First and Second
Lyapunov Theorems:

\noindent{\em the existence of a local Lyapunov function implies the
a.s. (open loop) stabilizability of (CSDE);
a strict Lyapunov function implies the
a.s. (open loop) asymptotic stabilizability}.

\noindent The same proof provides their global versions as well: if $V$
satisfies \eqref{cseq1.0}
in $\R^{N}\setminus\{0\}$ then  $(CSDE)$ is also a.s. (open loop)
{\em Lagrange stabilizable}, i.e., for all initial points $x$ there is a
control such that \eqref{cappa0} holds, and if $V$ is strict then the
system is {\em globally} a.s. (open loop) asymptotically stabilizable.
We also give sufficient conditions
for the stability of viable (controlled invariant) sets more
general than an equilibrium point, and for the a.s.  exponential stability.

These facts are much easier to prove if the Lyapunov function is
smooth, but this assumption is not necessary and would limit
considerably their
applicability.  The nonexistence of smooth Lyapunov functions is well
known in the deterministic case, see
\cite{kra, bacr} for stable uncontrolled systems and the surveys
\cite{son, bacr} for asymptotically stable controlled systems.  Here
we give an example of an uncontrolled degenerate diffusion process
that is a.s.  stable but cannot have a continuous Lyapunov function
(Example \ref{kras} in Section 6).  Moreover, in a companion paper
\cite{ces1} the second author proves a {\em converse Lyapunov
theorem}, stating that any a.s.  stabilizable system $(CSDE)$ has a
l.s.c.  local viscosity Lyapunov function.

All the results listed above are about open loop a.s.  stabilizability.
They raise the question of the existence of a stabilizing feedback.
Here we give an answer only for affine systems with a smooth strict Lyapunov
function.  We adapt Sontag's method \cite{son0} to the stochastic setting
and find an explicit formula for a feedback that renders the system
almost surely asymptotically stable.  The   feedback
stabilizability of controlled diffusions in case of nonsmooth Lyapunov
functions seems considerably harder and we are not aware of any paper
on the subject.

In the last section we study some simple applications and examples.
For instance, we consider a deterministic, asymptotically
controllable system $\dot{X}_t=f(X_t, \alpha_t)$ with Lyapunov pair
$(V,L)$ and look for conditions on a stochastic perturbation that keep
the system a.s.  stabilizable with the same Lyapunov function $V$ for
some $l\leq L$.

Our proof of the first Lyapunov-type theorem is based on the
observation that the infinitesimal decrease condition \eqref{cseq1.0}
has the rescaling property of the geometric PDEs arising in the level
set approach to front propagation (see, e.g., \cite{bcess, st2} and
the references therein), and on a recent result of the first author
and Jensen \cite{bj} on the viability, or controlled invariance, of
general closed sets for controlled diffusions (see \cite{au, ad3} and
the references therein for earlier work on viability for stochastic
processes).
For the second Lyapunov-type theorem we use also
martingale inequalities and other properties of diffusions.

The first Lyapunov-type theorem on local a.s.  stabilizability was
announced in \cite{bc} where we presented the simpler proof for
uncontrolled processes.  In the forthcoming paper \cite{ces2} the
second author shows that the existence of a continuous viscosity
solution of the Hamilton-Jacobi-Bellman inequality
$$
\max_{\alpha\in A} \left\{-DV(x)\cdot
f(x,\alpha)-trace\left[a(x,\alpha) D^2V(x)\right] \right\} \geq l(x) ,
$$
implies the open loop stabilizability in probability of $(CSDE)$.
Converse theorems in this setting are under investigation.

We conclude with some additional references.  Nonsmooth Lyapunov functions
for uncontrolled diffusion processes were studied by Ladde and Lakshmikantham
\cite{ll} with Dini-type derivatives along sample paths, and by Aubin
and Da Prato \cite{ad4} by means of a stochastic contingent epiderivative.
Recently Arnold and Schmalfuss \cite{as} gave an extension of
Lyapunov's Second Method to random dynamical systems.
Turning to deterministic controlled systems, we recall that Soravia
\cite{sor1} gave direct and inverse Lyapunov theorems for the open
loop stabilizability by means of viscosity solutions (in the more
general context of differential games), Sontag and Sussmann \cite{son0,
sosu} did it for the asymptotic controllability (i.e., asymptotic open
loop stabilizability) by using Dini directional derivatives.  Viscosity
methods for stability problems were also used in \cite{ksor, sor2,
gru}.  There is a large literature on feedback stabilization, see
\cite{art, son1, clss}, the surveys \cite{son, clsw, bacr}, and the
references therein.  We refer to \cite{cil, bcd} for the basic
theory of viscosity solutions, and to \cite{pll1, bcess, FSo, yz} for
its applications to deterministic and stochastic optimal control.


The paper is organized as follows.  In Section 2 we give the main
definitions and state the first and second Lyapunov-type theorems.
Section 3 recalls some viability theory and then gives the proofs of
the two main theorems.  Section 4 is on feedback stabilization of
affine systems with smooth Lyapunov functions.  Section 5 contains
some extensions to exponential stability, general equilibrium sets,
and target problems.  Section 7 is devoted to the examples.


\section{Lyapunov functions for a.s.  stabilizability and asymptotic
stabilizability}
We consider a controlled  Ito stochastic differential equation:
\[(CSDE)\left\{\begin{array}{l}
dX_t=f(X_t,\alpha_t)dt
+\sigma(X_t,\alpha_t)dB_t,\;\; t>0,\\ X_0=x,
\end{array}
\right. \]
where $B_t$ is an $M$-dimensional Brownian motion. Throughout the paper we
assume that
$f, \sigma$ are continuous functions defined in
$\R^N\times A$, where $A$ is a compact metric
 space, which take values, respectively, in $\R^N$ and in the space of
 $N\times M$ matrices, and satisfying
\begin{equation}\label{condition1}
|f(x,\alpha)-f(y,\alpha)|+
\Vert\sigma(x,\alpha)-\sigma(y,\alpha)\Vert\leq C|x-y|, \quad \forall
x,y\in\R^N, \quad \forall \alpha\in A.
%
%
\end{equation}
We adopt the definition of admissible control function,  or
admissible system, of Haussmann and Lepeltier \cite{hl} (Definition
2.2, page 853).
For a given $x\in \R^N $ we denote with ${\cal A}_x$ the set of
admissible control functions,
with
$\alpha_\cdot$ its generic element (although it is not a
standard function $\R\to A$), and with $X_\cdot$ the corresponding solution of
$(CSDE)$.
%
%
\noindent We define
$$a(x,\alpha):=\frac{1}{2}\sigma(x,\alpha)\sigma(x,\alpha)^T$$
and assume
\begin{equation}\label{convex}
\left\{(a(x,\alpha),f(x,\alpha)) \; : \; \alpha\in A \right\}\quad
\text{is convex for all } x\in \R^N .
\end{equation}
%
\begin{defi}[almost sure stabilizability]\label{def:as.stab}
The system $(CSDE)$ is \emph{almost surely (open-loop Lyapunov) stabilizable}
at the origin if for every
$\eta>0$ there exists $\delta>0$ such that, for any initial point $x$
with $|x|\leq \delta$, there exists an admissible control function
$\overline{\alpha}_{\cdot}\in{\cal A}_x$ whose corresponding
trajectory $\overline{X}_{\cdot}$ verifies $|\overline{X}_t|\leq \eta$
for all $t\geq 0$ almost surely.

The system is \emph{almost surely (open-loop) Lagrange stabilizable}, or it
has the property of \emph{uniform boundedness of trajectories}, if for
each $R>0$ there is $S>0$ such that for any initial point $x$ with $|x|\leq R$, there exists an
admissible control function $\overline{\alpha}_{\cdot}\in{\cal A}_x$
whose corresponding trajectory $\overline{X}_{\cdot}$
verifies $|\overline{X}_t|\leq S$ for all $t\geq 0$ almost surely.
\end{defi}
\begin{rem}\upshape
The a.s.  stabilizability implies that the origin is a {\em controlled
equilibrium} of $(CSDE)$, i.e., 
$$
\exists \,\oa\in A \,:\; f(0,\oa)=0, \; \sigma(0,\oa)=0.
$$
In fact, the definition
%
gives for any $\eps>0$
an admissible control such that the corresponding trajectory
starting at the origin satisfies a. s.
$|X_t|\leq \eps$ for all $t$,
so $\mathbf{E}_x\int_0^{+\infty} |X_t|e^{-\lambda t}dt \leq
\frac{\eps}{\lambda}$ for any $\lambda>0$.
Then
$\inf_{\alpha_.\in {\cal A}_x} \mathbf{E}_x\int_0^{+\infty}
|X_t|
e^{-\lambda t}dt=0
$.
The convexity assumption (\ref{convex}) and an existence theorem
for optimal controls \cite{hl} imply that the $\inf$ is attained,
and the minimizing control produces a trajectory satisfying a.s.
$|X_t|=0$ for all $t\geq 0$.
 The conclusion follows from standard
properties of stochastic differential equations.
\end{rem}
\begin{rem}\upshape
As it is common in the modern deterministic stability theory, the
previous definitions can be reformulated in terms of the {\em
comparison functions} introduced by Hahn \cite{ha}.  We will use
the class $\mathcal{K}$ of continuous functions
$\gamma:[0,+\infty)\rightarrow[0,+\infty)$ strictly increasing and
such that $\gamma(0)=0$, and the class $\mathcal{K}_{\infty}$ of
functions $\gamma\in\mathcal{K}$ such that $\lim_{r\rightarrow
+\infty}\gamma(r)=+\infty$.
%
The system $(CSDE)$ is a.s. (open-loop) stabilizable at 0 if  there exists
$\gamma\in\mathcal{K}$ and $\delta_o>0$ such that for any starting point
$x$ with $|x|\leq \delta_o$
\begin{equation}
\label{cappa}
\exists \,
\overline{\alpha}_{\cdot}\in{\cal A}_x \,: \quad  |\overline{X}_t|\leq
\gamma(|x|)\quad \forall t\geq 0 \,\text{ a.s.,}
\end {equation}
 where $\overline{X}_t$ is the trajectory corresponding to
$\overline{\alpha}_{\cdot}$.  If (\ref{cappa}) holds for some
$\gamma\in \mathcal{K}_{\infty}$ and for all $x\in\R^N$, then the
system is also a.s.  (open-loop) Lagrange stabilizable.
\end{rem}
\begin{defi}[a.s. asymptotic stabilizability]
The system $(CSDE)$ is \emph{almost surely (open loop) locally
asymptotically stabilizable} (or \emph{a.  s.  locally
asymptotically controllable}) at the origin if for every $\eta>0$
there exists $\delta>0$ such that, for all $|x|\leq \delta$, there
exists an admissible control function
$\overline{\alpha}_{\cdot}\in{\cal A}_{x}$ whose corresponding
trajectory $\overline{X}_{\cdot}$ %
%
verifies almost surely
$$
|\overline{X}_t|\leq \eta\quad \forall t\geq 0,
\qquad\lim_{t\rightarrow+\infty}|\overline{X}_t|=0.
$$

The system is \emph{a.  s. (open loop) globally asymptotically
stabilizable} (or {\em a. s. asymptotically controllable})  at the
origin if there is $\gamma\in \mathcal{K}_{\infty}$ and for all
$x\in\R^N$ there exists $\overline{\alpha}_{\cdot}\in{\cal A}_{x}$
whose trajectory $\overline{X}_{\cdot}$ satisfies almost surely $$
|\overline{X}_t|\leq \gamma(|x|) \quad \forall t\geq 0,
\qquad\lim_{t\rightarrow +\infty}|\overline{X}_t|=0. $$

\end{defi}

Next we give the appropriate definition of Lyapunov function
for the study of almost sure stabilizability.
We recall the definition of second order semijet of a l.s.c.  function
$V$ at a point $x$:
\begin{eqnarray*}
{\cal J}^{2,-} V(x) := \left\{ (p,Y)\in\R^N\times S(N): \text{ for
}y\rightarrow x \right.  \\
 \left.
V(y)\geq V(x)+p\cdot (y-x)+\frac{1}{2} (y-x)\cdot
Y(y-x)+o(|y-x|^2)\right\}.
\end{eqnarray*}

\begin{defi}[control Lyapunov function]\label{cliap}
Let ${\cal O}\subseteq \R^N$ be an open set containing the origin.
A function $V : {\cal O}\rightarrow[0,+\infty)$ is a \emph{control
Lyapunov function} for the a.s.  stability of $(CSDE)$ if

\noindent (i) $V$ is lower semicontinuous,

\noindent(ii) $V$ is  continuous at $0$ and {\em positive
definite}, i.e., $V(0)=0$ and $V(x)>0$ \\for all $x\neq 0$,

\noindent(iii) $V$ is {\em proper}, i.e.,
$\lim_{|x|\rightarrow +\infty}V(x)=+\infty$ or, equivalently, the
level sets \\ $\left\{x| V(x)\leq\mu\right\}$ are bounded for every $\mu\in
[0,\infty)$;


\noindent (iv) for all $x\in{\cal O}\setminus\{0\}$ and
$(p,Y)\in{\cal J}^{2,-}V(x)$ there exists $\overline{\alpha}\in A$
such that
\begin{equation}
\label{ceq} \sigma(x,\overline{\alpha})^Tp=0 \;\ \ \  \,\mbox{ and
}\,\ \ \ -p\cdot
f(x,\overline{\alpha})-trace\left[a(x,\overline{\alpha}) Y\right]
\geq 0 .
\end{equation}
%
\end{defi}
\begin{rem}\upshape
The conditions $(ii)$ and $(iii)$ in the previous definition can be stated as
\begin{equation}
\label{cappainfty}
 \exists\,\gamma_1,\gamma_2 \in\mathcal{K}_{\infty}\;:\;
\gamma_1(|x|)\leq V(x) \leq \gamma_2(|x|)\; \forall\, x\in\R^N.
\end {equation}
 %
Therefore the level sets $\{V(x)\leq\mu \}$ of the Lyapunov function form a
basis of neighborhoods of $0$.
\end{rem}
\begin{rem}\upshape
If the dispersion matrix $\sigma$ does not depend on the control,
then condition $(iv)$ can be reformulated as follows:

\noindent $V$ is a solution in viscosity sense in ${\cal
O}\setminus\{0\}$ of the system:
\[\left\{\begin{array}{l}
\sigma(x)^TDV(x)=0 \\ \max_{\alpha\in A}\left\{-DV(x)\cdot
f(x,\alpha)-trace\left[a(x,\alpha) D^2V(x)\right] \right\}\geq 0
\end{array}\right.\]

In the general case, we can observe that if the condition $(iv)$
holds, then in particular $V$ is a viscosity supersolution of \begin{equation}
\label{ceq0}
\max_{\alpha\in A}\left\{-DV(x)\cdot
f(x,\alpha)-trace\left[a(x,\alpha) D^2V(x)\right] \right\}= 0\end{equation}
Moreover, if the function $V$ is at least differentiable, then the
condition $(iv)$ can be stated more concisely as:

\noindent $V$ is a supersolution in viscosity sense in ${\cal
O}\setminus\{0\}$ of the equation:
\[
\mbox{$\max_{\left\{\alpha\in A\ |\ \sigma(x,\alpha)^TDV(x)=0\right\}}$}
\left\{-DV(x)\cdot f(x,\alpha)-trace\left[a(x,\alpha)
D^2V(x)\right] \right\}= 0\]


\end{rem}
\begin{defi}[strict control Lyapunov function]\label{scliap}
A function $V:{\cal O}\rightarrow [0,+\infty)$ is a \emph{strict
control Lyapunov function} for the a.s.  stability of $(CSDE)$ if
it satisfies conditions (i), (ii), (iii) in the Definition
\ref{cliap} and


\noindent (iv)$^{\prime}$ for all $x\in{\cal O}\setminus\{0\}$ and
$(p,Y)\in{\cal J}^{2,-}V(x)$ there exists $\overline{\alpha}\in A$
such that
\begin{equation}
\label{cseq} \sigma^T(x,\overline{\alpha})p=0 \  \mbox{ and }\ \
-p\cdot f(x,\overline{\alpha})-trace\left[a(x,\overline{\alpha})
Y\right] - l(x)\geq 0
\end{equation}
for some positive definite and Lipschitz continuous $l:{\cal O}\to \R$.
\end{defi}
\begin{rem}\upshape
In the inequality in $(iv)^{\prime}$ we could take

$$ p\cdot f(x,\oa)-trace\left[a(x,\oa) Y\right]
-l(x,\oa)\geq 0 
$$
for some continuous $l:{\cal O}\times A\to \R$, Lipschitz
continuous in $x$ uniformly in $\alpha$, with $l(x,A)$ convex for
all $x\in {\cal O}$, and such that $\lt(x):=\min_{\alpha\in
A}l(x,\alpha)$ is positive definite. However, this would not
increase the generality of the definition, because $V$ would also
satisfy condition (\ref{cseq}) with $l$ replaced by $\lt$.
\end{rem}

\noindent
Our main results are the following versions for stochastic controlled
systems of the First and the Second Lyapunov Theorem.
\begin{teo}
    [a.s.  stabilizability] \label{thmliap1}
Assume (\ref{condition1}), (\ref{convex}), and the existence of a control
Lyapunov function $V$.  Then

\noindent (i) the system $(CSDE)$ is almost surely stabilizable at the origin;

\noindent (ii) if, in addition, the domain ${\cal O}$ of $V$ is all $\R^N$,
the system is also a.s.  Lagrange stabilizable and for all $x\in\R^N$ there
exists $\oa_.\in {\cal A}_x$ such that the corresponding trajectory
$\ox_.$ satisfies
\begin{equation}
\label{stimaglob}
\vert \ox_t \vert \leq \gamma_1^{-1}(\gamma_2(|x|))\quad \forall\,
t\geq 0\quad{a.s.}
\end{equation}
with $\gamma_1, \gamma_2\in {\cal K}_{\infty}$ verifying (\ref{cappainfty}).
\end{teo}
\begin{teo}
    [a.s. asymptotic  stabilizability]\label{thmliap2}
Assume (\ref{condition1}), (\ref{convex}), and the existence of a
strict control Lyapunov function $V$.  Then

\noindent (i) the system $(CSDE)$ is a.s. locally asymptotically stabilizable at the
origin;

\noindent (ii) if, in addition, the domain ${\cal O}$ of $V$ is all $\R^N$,
the system is a.s. globally asymptotically stabilizable.
\end{teo}

\section{A viability theorem and the proofs of stabilizability}

In this section we prove the Theorems \ref{thmliap1} and
\ref{thmliap2}.  Our main tool is a recent result in \cite{bj} about
the almost sure viability (named also {\em controlled invariance} and
{\em weak invariance})
 of an arbitrary closed set for a controlled diffusion process.
 (See \cite{au, ad3} and the references therein for earlier related results).
\begin{defi}[viable set] A closed set $K\subset\R^N$ is {\em viable}
or {\em controlled invariant} or {\em weakly invariant}
for the stochastic system $(CSDE)$ if for
all initial points $x\in K$ there exists an admissible control
$\alpha_.\in{\cal A}_x$ such that the corresponding trajectory $X_.$
satisfies $X_t\in K$ for all $t>0$ almost surely.
\end{defi}
\noindent It is easy to see from its very definition that the a.s.
stabilizability follows from the viability of all the sublevel
sets of any function satisfying the conditions $(i)-(iii)$ of
Definition \ref{cliap}.  The next result gives a geometric
characterization of viable sets.  It will allow us to check that
the sublevel sets of a control Lyapunov function are viable by
means of the condition $(iv)$ in Definition \ref{cliap}. The
Nagumo-type geometric condition in the viability theorem is given
in terms of the following {\em second order normal cone} to a
closed set $K\subset\R^N$, first introduced in \cite{bg},
\begin{eqnarray*}
{\cal N}^2_K (x)\!\! & := &\!\! \{(p,Y)\in I\!\!R^N \times S(N):
\;\;{\rm for}\;\; y \rightarrow x , \;\;y\in K ,\;\; \\
&    & p\cdot
(y-x)+\frac{1}{2} (y-x)\cdot Y(y-x)\geq o(|y-x|^2) \;\;\}
\end{eqnarray*}
where $S(N)$ is the set of symmetric $N\times N$ matrices.  Note that,
if $(p,Y)\in{\cal N}^2_K (x)$ and $x\in\partial K$, the vector $p$ is a
generalized (proximal or Bony) interior normal to the set $K$ at $x$.
In particular, if $\partial K$ is a smooth surface in a neighborhood
of $x$, $p/|p|$ is the interior normal and $Y$ is related to the second
fundamental form of $\partial K$ at $x$, see \cite{bg}.

 %
\begin{teo}[Viability theorem \cite{bj}]
\label{cor.convex} Assume (\ref{condition1}) and (\ref{convex}).  Then a
closed set $K\subseteq\R^N$
%
%
is viable for $(CSDE)$ if and only if
\begin{equation}
\forall x\in\partial K, \; \forall (p,Y)\in {\cal N}^2_K (x),\;
\exists \alpha\in A\;:\;\; f(x,\alpha)\cdot p +
trace\left[a(x,\alpha) Y\right]\geq 0 .\label{viabcond}
\end{equation} 
\end{teo}

%

\noindent The second tool for the proof of the Lyapunov-type
Theorem \ref{thmliap1} is the following lemma on the change of
unknown for second order partial differential equations.  It says
that the Hamilton-Ja\-co\-bi-Bell\-man inequality in condition
(\ref{ceq}) in the definition of control Lyapunov function behaves
as a \emph{geometric equation} if the unknown satisfies also the
condition in (\ref{ceq})
of orthogonality between its gradient and the columns of
the dispersion matrix $\sigma$.
We refer the interested reader to the chapters by Evans and by
Souganidis in the book \cite{bcess} for an introduction to the
geometric PDEs of the theory of front propagation.
\begin{lemma}\label{lemma2}
Let $v$ satisfy condition (\ref{ceq}) for all $(p,Y)\in{\cal
J}^{2,-}V(x)$, $x\in \R^N\setminus\{0\}$.  Let $\phi$ be a twice
continuously differentiable strictly increasing real map.  Then
$w=\phi\circ v$ is a viscosity supersolution of
\begin{equation}\label{lemmaeq} \max_{\alpha\in
A}\left\{-DV(x)\cdot f(x,\alpha)-trace\left[a(x,\alpha)
D^2V(x)\right] \right\}= 0\end{equation}
\end{lemma}
\begin{proof} It is easy to check that, if $(p,Y)\in{\cal J}^{2,-}w(x)$, then
\[
\left(\psi'(w(x))p , \psi'(w(x))Y+\psi''(w(x))p\otimes p\right) \in {\cal
J}^{2,-}v(x),
\]
where $\psi$ is the inverse of $\phi$ and $p\otimes p$ is the
$N\times N$ matrix whose $(i,j)$ entry is $p_ip_j$.  Then, for
$(p,Y)\in{\cal J}^{2,-}w(x)$ and $x\ne 0$ there exists $\oa$ such
that
\[
\left\{ -\psi'(w(x))p\cdot
f(x,\oa)-trace\left[a(x,\oa)\cdot\left(\psi'(w(x))Y
+\psi''(w(x))p\otimes p \right)\right] \right\}\geq 0.
\]
and
\[
trace \left[a(x,\overline{\alpha})\cdot\psi''(w(x))p\otimes p
\right]= \frac{\psi''(w(x))}{(\psi'(w(x)))^2}
|\sigma(x,\overline{\alpha})^T\psi'(w(x))p|^2=0 .
\]
Therefore
\[
-\psi'(w(x))p\cdot f(x,\overline{\alpha})-trace
\left[a(x,\overline{\alpha}) \cdot\psi'(w(x))Y \right]\geq 0
\]
and we can conclude that
\[
\sup_{\alpha\in A} \left\{ -p\cdot
f(x,\alpha)-trace\left[a(x,\alpha)\cdot Y \right] \right\}\geq 0 .
\]
\end{proof}

\begin{proof} {\it of Theorem \ref{thmliap1}.}
We begin with the proof of $(ii)$. We fix an arbitrary $\mu>0$ and consider the
sublevel set of the function $V$
\[
K := \{x \,|\, V(x)\leq\mu\}.
\]
We claim that $K$ is viable.  Then for all initial points $x\in\R^N$
there exists $\oa_.\in{\cal A}_x$ such that the associated trajectory
$\ox_.$ satisfies
$$
\gamma_1(\vert \ox_t\vert)\leq V(\ox_t) \leq V(x) \leq
\gamma_1(|x|),\quad \forall t\geq 0\quad\text{a.s.},
$$
which gives the estimate (\ref{stimaglob}).  Then the system is a.s.
stabilizable and Lagrange stabilizable because $\gamma_1^{-1}\circ
\gamma_2\in{\cal K}_{\infty}$.

To prove
that $K$ is viable we will check the condition (\ref{viabcond}) of the
Viability Theorem \ref{cor.convex}.  %
For a given $\lambda>0$ we define the nondecreasing continuous real function
\[
\psi_{\lambda}(t)= \left\{\begin{array}{ll} 0 , & t\leq \mu ,
\\ \lambda (t-\mu ) , & \mu\leq t ,
\leq \mu+\frac{1}{\lambda} , \\ 1 , & t\geq \mu+\frac{1}{\lambda} .
\end{array} \right.
\]
We claim that the  function $\psi_{\lambda}\circ V$ is a viscosity
supersolution of equation (\ref{lemmaeq}) for every $\lambda$.  To
prove the claim we choose a sequence $\psi_n$ of strictly
increasing, smooth real maps that converge uniformly on compact
sets to $\psi_{\lambda}$. Then, for every $n$, the map
$\psi_n\circ V$ is a viscosity supersolution of equation
(\ref{lemmaeq}) by Lemma \ref{lemma2}.  By the stability of
viscosity supersolutions with respect to uniform convergence
we get the claim.

Next we observe that the net $\psi_{\lambda}\circ V$ is increasing and
converges as $\lambda\rightarrow +\infty$ to the indicator function
\[
C(x)=\left\{\begin{array}{ll} 0 , & x\in K , \\ 1 , & x\not\in K .
\end{array} \right.
\]
Viscosity supersolutions are stable with respect to the pointwise
increasing convergence (see, e.g., Prop.  V.2.16, p.  306, of
\cite{bcd}). Therefore the indicator function $C$ of $K$ is a
viscosity supersolution of equation (\ref{lemmaeq}).
 From the definitions it is easy to check that
\[
{\cal J}^{2,-} C(x)=-{\cal N}^{2}_K (x),\quad \forall \, x\in \partial
K .
\]
By plugging this formula into the equation (\ref{lemmaeq}) we
obtain exactly the condition (\ref{viabcond}) of the Viability
Theorem and complete the proof of $(ii)$.

To prove $(i)$ we choose $\overline{\mu}>0$ small enough so that
$K:=\{x\in {\cal O}\,:\, V(x)\leq
\mu\}$, for
 $\mu\leq\overline{\mu}$, is closed in $\R^N$ (for instance,
$\overline{\mu}<\inf_{y\in\partial\mathcal{O}}\liminf_{x\to y} V(x)$).
Then the
preceding part of this proof gives the viability of $K$ and the
estimate (\ref{stimaglob}) for all $x$ such that
$V(x)\leq\overline{\mu}$.  Therefore, for some $\delta_o>0$,
(\ref{stimaglob}) holds for all $x$ with $|x|\leq\delta_o$, and this
gives the a.s.  stabilizability of the origin.
\end{proof}

Next we give the proof of Theorem \ref{thmliap2} about
asymptotic stability.  It is obtained by first applying the previous Theorem
\ref{thmliap1} to a new system with an extra variable, and then using
martingale inequalities as, e.g., in \cite{deng}.

%
\begin{proof} {\it of Theorem \ref{thmliap2}.}
We consider the differential system
\[
\left\{
\begin{array}{cclcc}
dX_t & = & f(X_t,\alpha_t)dt & +&\sigma(X_t,\alpha_t)dB_t\\ dZ_t &
= & l(X_t)dt& &
\end{array} \right.
\]
with initial data $X_0=x$ and $Z_0=0$.  We rewrite this system in
$\R^{N+1}$ as
\[
(CSDE2)\left\{\begin{array}{l}
d(X_t,Z_t)=\overline{f}(X_t,Z_t,\alpha_t)dt+
\overline{\sigma}(X_t,Z_t,\alpha_t)d(B_t,0),\;\; t>0,\\
(X_0,Z_0)=(x,0).
\end{array}
\right.
\]
where
$\overline{f}(x,z,\alpha)=(f(x,\alpha),l(x))$ and
$\overline{\sigma}(x,z,\alpha)= (\sigma(x,\alpha),0)$.
Clearly it
satisfies the conditions (\ref{condition1}) and  (\ref{convex}).
Let us consider the function
\begin{eqnarray*}
 W(x,z):&
{\cal O}\times \R\rightarrow &\R\\
& (x,z)\longmapsto & V(x)+|z|
\end{eqnarray*}
We claim that it is a Lyapunov function for $(CSDE2)$.  In fact,
$W$ is positive definite (because $W\geq 0$ and $W=0$ only for
$(x,z)=(0,0)$); $W$ is lower semicontinuous, continuous at $(0,0)$
and proper since $V$ is so.  We have only to prove that $W$
satisfies condition (\ref{ceq}). Fix $x\ne 0$ and $(x,z)$ with
$z>0$ and a smooth function $\phi$ such that $W-\phi$ has a local
minimum at $(x,z)$, i.e.,
\[
V(x)+z-\phi(x,z)\leq V(y)+w-\phi(y,w)
\]
for every $(y,w)$, $w>0$ in a neighborhood of $(x,z)$.
If we choose $w=z$
we get a minimum in $x$ for the function $V(\cdot)-\phi(\cdot,z)$,
 therefore $(D_x\phi(x,z),D^2_{xx}\phi(x,z))\in
J^{2,-}V(x)$; if we choose $y=x$ we find a minimum in $z$ for the
smooth function $w\longmapsto w-\phi(x,w)$, so $D_z\phi(x,z)=1$.
Then there exists $\overline{\alpha}\in A$ such that
$(\sigma(x,\overline{\alpha}),0)^T (D_x\phi(x,z),1) =0$ and
$$\left\{-D\phi(x,z)\cdot \overline{f}(x,z,\oa)-trace
\left[\overline{a}(x,z,\oa)D^2\phi(x,z)\right]\right\}= $$ $$=
\left\{-(D_x\phi(x,z),1)\left(\begin{array}{c} f(x,\oa)\\ l(x)
\end{array} \right)-trace \left[ \left(
\begin{array}{cc} a(x,\oa) & 0\\ 0& 0 \end{array}
\right)D^2\phi(x,z)  \right]\right\} = $$ $$=
\left\{-D_x\phi(x,z)\cdot f(x,\oa)-trace \left[
a(x,\oa)D^2_{xx}\phi(x,z)\right]\right\}-l(x)\geq 0 , $$
%
since $V$ is a strict Lyapunov function. Now fix $(x,z)$ with $z<
0$ and let $\phi$ be a smooth function such that
%
\[
V(x)-z-\phi(x,z)\leq V(y)-w-\phi(y,w)
\]
for every $(y,w)$, $w<0$ in a neighborhood of
$(x,z)$.
We argue as before and now get that there exists
$\overline{\alpha}\in A$ such that
$
(\sigma(x,\overline{\alpha}),0)^T\cdot (D_x\phi(x,z),-1) =0
$
and $$ -D_x\phi(x,z)\cdot
f(x,\overline{\alpha})-trace\left[a(x,\overline{\alpha})
D^2_{xx}\phi(x,z)\right] + l(x) >
 $$
 $$
 >-D_x\phi(x,z)\cdot
 f(x,\overline{\alpha})-trace\left[a(x,\overline{\alpha})D^2_{xx}\phi(x,z)
 \right]-l(x) \geq 0 .
 $$
because $l$ is positive and $V$ is a Lyapunov function.

Finally, we consider $(x,0)$ and a smooth function  $\phi$ such
that
\[
V(x)-\phi(x,0)\leq V(y)-w-\phi(y,w)
\]
for every $(y,w)$, $w<0$ in a neighborhood of
$(x,0)$ and
\[
V(x)-\phi(x,0)\leq V(y)+w-\phi(y,w)
\]
for all $(y,w)$, $w>0$ in a neighborhood of $(x,0)$.
Then $(D_x\phi(x,z),D^2_{xx}\phi(x,z))\in J^{2,-}V(x)$,
$D_z\phi(x,0)
\geq -1$, and
$D_z\phi(x,0)
\leq 1$.  Therefore there exists $\overline{\alpha}\in A$ such
that
$
(\sigma(x,\overline{\alpha}),0)^T\cdot (D_x\phi(x,z),D_z\phi(x,z))
=0  $
and $$\left\{-D\phi(x,z)\cdot \overline{f}(x,z,\oa)-trace
\left[\overline{a}(x,z,\oa)D^2\phi(x,z)\right]\right\}= $$
$$ =\left\{-D_x\phi(x,z)\cdot f(x,\oa)-trace \left[
a(x,\oa)D^2_{xx}\phi(x,z)\right]\right\}-D_z\phi(x,z)l(x)\geq $$
$$ =\left\{-D_x\phi(x,z)\cdot f(x,\oa)-trace \left[
a(x,\oa)D^2_{xx}\phi(x,z)\right]\right\}-l(x)\geq 0 , $$

This completes the proof of the claim,  so we can apply Theorem
\ref{thmliap1} to get for every $x\in {\cal O}$  an admissible
control $\oa_.\in{\cal A}_x$ such that the corresponding
trajectory $(\ox_.,\overline{Z}_.)$ of (CSDE2) with initial data
$(x,0)$ remains almost surely in the level set $K=\{(y,w)\in {\cal
O}\times\R \,|\, W(y,w)\leq W(x,0)\}$.  Then, for all $t\geq 0$
and almost surely, $\ox_t\in {\cal O}$,
\[
W(\ox_t,\overline{Z}_t)=V(\ox_t) + \overline{Z}_t = V(\ox_t)+
\int_0^t\!\!l(\ox_s)ds\leq W(x,0)=V(x) ,
\]
and
\begin{equation}\label{asineq}
0\leq V(\ox_t)\leq V(x)-\int_0^t\!\!l(\ox_s)ds .  
\end{equation}
In particular, since $l\geq 0$, for some $r>0$, $|\ox_t|\leq r$ for all
$t$ almost surely.


Next we claim that $l(\ox_t)\rightarrow 0$ almost surely as
$t\rightarrow +\infty$.  Let us assume by contradiction that the claim
is not true: then there exist $\eps>0$, a subset $\Omega_\eps\subseteq
\Omega$ with $\mathbf{P}(\Omega_\eps)>0$, and for every $\omega\in
\Omega_\eps$ 
a sequence $t_n(\omega)\rightarrow +\infty$ such that
$l(\ox_{t_n}(\omega)) >\eps$.  We define
 $$
 F(r):=\max_{|x|\leq r, \alpha\in A}|f(x,\alpha)| ,
 \qquad\Sigma(r)=\max_{|x|\leq r, \alpha\in A} \|\sigma(x,\alpha)\| .
$$
 We compute
\[
\mathbf{E}\left\{\sup_{t\leq s\leq t+h}|\ox_s-\ox_t|^2
\right\}=
\]
\[
=\mathbf{E}\left\{\sup_{t\leq s\leq t+h}|\int_t^s
f(\ox_u,\oa_u)du+\int_t^s \sigma(\ox_u,\oa_u)dB_u|^2 \right\}\leq
\]
\[
\leq  2\mathbf{E}\left\{\sup_{t\leq s\leq t+h}|\int_t^s
f(\ox_u,\oa_u)du|^2\right\}+2 \mathbf{E}\left\{\sup_{t\leq s\leq
t+h}|\int_t^s \sigma(\ox_u,\oa_u)dB_u|^2 \right\}\leq
\]
\[
\leq 2 F^2(r)h^2+2 \mathbf{E}\left\{\sup_{t\leq s\leq t+h}|\int_t^s
\sigma(\ox_u,\oa_u)dB_u|^2\right\} =: K .
\]
By Theorem 3.4 in \cite{doob} (the process $|\int_t^s
\sigma(\ox_u,\oa_u)dB_u|$ is a positive semimartingale) we get
\[
K \leq
2 F^2(r)h^2+8 \sup_{t\leq s\leq t+h}\mathbf{E}\left\{|\int_t^{s}
\sigma(\ox_u,\oa_u)dB_u|^2\right\}
\]
and by the Ito isometry
\[
K \leq 2 F^2(r)h^2+8 \mathbf{E}\left\{\int_t^{t+h}
|\sigma(\ox_u,\oa_u)|^2du\right\}\leq 2 F^2(r)h^2+8 \Sigma^2(r)h .
\]
Then, Chebyshev inequality gives
$$
\mathbf{P}\left\{ \sup_{t\leq s\leq t+h}
|\ox_s-\ox_t|>k \right\} \leq \frac{\mathbf{E} \left\{ \sup_{t\leq
s\leq t+h} |\ox_s-\ox_t|^2\right\}}{k^2} \leq
$$
\[
\leq \frac{2 F^2(r)h^2+8 \Sigma^2(r)h}{k^2} .
\]
Since $l$ is continuous, we can fix $\delta$ such that
$|l(x)-l(y)|\leq\frac{\eps}{2}$ if $|x-y|\leq\delta$ and
$|x|, |y| \leq r$.  We define
$$
C := \left\{\omega\in\Omega \;:\; \sup_{0\leq s\leq h} |\ox_s-x|
\leq\delta\right\}
$$
and choose $0<k<\mathbf{P}(\Omega_\eps)$
and $h>0$ depending on $\delta$ and $\eps$ such that
\[
\mathbf{P}_x(C) \geq 1- \frac{2
F^2(r)h^2+8 \Sigma^2(r)h}{\delta^2}\geq 1+k-\mathbf{P}(\Omega_\eps).
\]
By the uniform continuity of $l$ then the set \[B:=
\left\{\omega\in\Omega \;:\; \sup_{0\leq s\leq h} |l(\ox_s)-l(x)|
\leq\eps/2\right\}\] contains $C$ and then \begin{equation}
\label{uni}\mathbf{P}_x(B)
\geq 1+k-\mathbf{P}(\Omega_\eps).\end{equation}
From the inequality (\ref{asineq}), letting $t\rightarrow \infty$,
we get
\[V(x)\geq
\mathbf{E}_x\int_0^{+\infty}\!\!l(\ox_s) ds \geq
\int_{\Omega_\eps} \int_0^{+\infty}\!\!l(\ox_s) \, ds\,
d\mathbf{P}\geq \int_{\Omega_\eps} \sum_{n}
\int_{t_n(\omega)} ^{t_{n}(\omega)+h}\!\!  l(\ox_s) ds\,
d\mathbf{P}\geq
\]\[\geq \int_{\Omega_\eps} \sum_{n}h \inf_{[t_n(\omega),t_n(\omega)+h]} l(\ox_t)\geq  h \sum_{n}\int_{\Omega_\eps}
\inf_{[t_n(\omega),t_n(\omega)+h]} l(\ox_t)d\mathbf{P}\geq  \]
\[\geq  h \sum_{n}\frac{\eps}{2} \mathbf{P}
\left[\left(\sup_{0\leq s\leq h}
|l(\ox_s)-l(x)|\leq\eps/2\ |\ x=\ox_{t_n}\right)\cap\Omega_\eps\right]
\]  by the strong Markov property
of the solutions of (CSDE) the estimate (\ref{uni}) gives
$\mathbf{P}\left(\sup_{0\leq s\leq h}
|l(\ox_s)-l(x)|\leq\eps/2\ |\ x=\ox_{t_n}\right)\geq
1+k-\mathbf{P}_x(\Omega_\eps)$ for every $n$. Therefore
$ \mathbf{P}
\left[\left(\sup_{0\leq s\leq h}
|l(\ox_s)-l(x)|\leq\eps/2\ |\ x=\ox_{t_n}\right)\cap\Omega_\eps\right]\geq k$
for every $n$: so, by the previous inequality, we get
\[V(x)\geq   h \sum_{n}\frac{\eps}{2}k=+\infty\]
This  gives a contradiction: then, for every $\eps>0$,
$\mathbf{P}(\Omega_\eps)=0$.
 We have proved that $l(\ox_t)\rightarrow 0$ almost surely as
$t\rightarrow +\infty$, now the positive definiteness of $l$ implies that
$|\ox_t|\rightarrow 0$ almost surely as $t\rightarrow +\infty$.
\end{proof}
\begin{rem}
\upshape
If the function $l$  is only nonnegative semidefinite
the proof of the last theorem gives, for any $x$, a control $\oa_.$ whose
trajectory $\ox_t$
satisfies a.s.  $V(\ox_t)\leq V(x)$ and $l(\ox_t)\rightarrow 0$ as
$t\rightarrow +\infty$.  Then the set ${\cal L}:=\{y\ | \ l(y)=0\}$ is
an attractor, for a suitable choice of the control, in the sense that
$dist(\ox_t,{\cal L}) \rightarrow 0$ a.s. as $t\rightarrow +\infty$.
For uncontrolled diffusion processes results of this kind can be found in
\cite{mao} and \cite{deng} and they are considered as stochastic versions of
a theorem by La Salle.  The earlier paper of Kushner \cite{kus2} studies
also a stochastic version of the La Salle invariance principle,
namely, that the omega limit set of the process is an invariant subset
of ${\cal L}$, in a suitable sense.
%
\end{rem}

\section{A.s.  feedback stabilization of affine systems}

In this section we give a result on the \emph{feedback
stabilizability} of systems affine in the control in the case there
exists a smooth strict control Lyapunov function.  It is an analogue for
the a.s.  stability of a celebrated theorem of Artstein \cite{art} and
Sontag \cite{son1} for deterministic systems, extended by Florchinger
\cite{flo} to the stability in probability of controlled diffusions.


We begin with the simple case of a single-input affine system with
uncontrolled diffusion, that is,
\begin{equation}\label{affine1}
 dX_t=\left(f(X_t)+\alpha_t g(X_t)\right)dt+
\sigma(X_t) dB_t ,
\end{equation}
where $f,g, \sigma$ are vector fields in $\R^N$ with $f(0)=0$ and
$\sigma(0)=0$, $B_t$ is a $1$-dimensional Brownian motion, and the
control $\alpha_t$ takes values in $\R$.
We seek a
function $k : \R^N\to \R$, at least continuous in $\R^N\setminus \{0\}$, such
that the origin is a.s.  asymptotically stable for the stochastic
differential equation
\begin{equation}\label{sde}
 dX_t=\left(f(X_t)+ k(X_t) g(X_t)\right)dt+
\sigma(X_t) dB_t .
\end{equation}
Then $k$ is called an \emph{a.s.  asymptotically stabilizing feedback} for
the control system (\ref{affine1}).

If there are no constraints on the control, a smooth strict control Lyapunov
function $V$ satisfies, in $\R^N\setminus\{0\}$,
$$
f\cdot DV + trace\left[\frac{1}{2}\sigma\sigma^T D^2V\right] +
\inf_{\alpha\in \R}\left\{\alpha g\cdot DV \right\}\leq - l , \qquad \sigma \cdot DV =
0 .
$$
Set $\gamma(x):= f\cdot DV + trace\left[\sigma\sigma^T D^2V\right]/2
+l/2$ and observe that the inequality for $V$ means
$$
g(x)\cdot DV(x) = 0 \quad \Rightarrow \quad \gamma(x)\leq -l(x)/2 < 0 .
$$
It is clear that $k(x):= -\gamma(x)/g(x)\cdot DV(x)$, $k(x):=0$ if
$g(x)\cdot DV(x) = 0$, could be a stabilizing feedback, but it is
discontinuous where $g(x)\cdot DV(x)$ vanishes.  If this case occurs
we build a continuous feedback by means of Sontag's universal formula
\cite{son1}, i.e.,
\begin{equation}
\label{formulafeed}
k(x) := - \frac{\gamma(x) +
\sqrt{\gamma^2(x) + (g(x)\cdot DV(x))^4}}{g(x)\cdot DV(x)}, \quad
\text{if } g(x)\cdot DV(x)\ne 0,
\end {equation}
 and $k(x)=0$ if $g(x)\cdot DV(x)=0$.  By the argument in
\cite{son1} $k\in C(\R^N\setminus\{0\})$ if $V\in
C^2(\R^N\setminus\{0\})$, and $k\in C^1(\R^N\setminus\{0\})$ if $f, g,
l$ are of class $C^1$ and $V\in C^3(\R^N\setminus\{0\})$.  Moreover
$$
(f + kg)\cdot DV + trace\left[\frac{1}{2}\sigma\sigma^T D^2V\right] \leq -
\frac{l}{2} , \qquad \sigma \cdot DV = 0
$$
in $\R^N\setminus \{0\}$, so $V$ is a strict Lyapunov function for
(\ref{sde}) and the origin is a.s.  asymptotically stable.  In
conclusion, $k$ is a stabilizing feedback for the affine control
system (\ref{affine1}).

If the control must satisfy a hard constraint, say $\alpha\in [-1,
1]$, it is not hard to check that $k(x)$ can be used in a neighborhood
of the origin provided that $DV$ and $D^2V$ are bounded near 0 and
either $g(x)\to 0$ or $DV(x)\to 0$ as $x\to 0$.

Next we use the same idea for the more general system with both the
drift and the diffusion terms affine in the control
\begin{equation}\label{affine2}
dX_t=\left(f(X_t)+\sum_{i=1}^{P-1} \alpha^i_t g_i(X_t)\right) dt+
\left(\sigma(X_t) + \alpha^P_t \tau(X_t) \right) dB_t .
\end{equation}
where $f,g_i, \sigma,
\tau$ are vector fields in $\R^N$, $B_t$ is a standard $1$-dimensional
Brownian motion, and the controls $\alpha^i_t$, $i=1,...,P$, are $\R$-valued.
The existence of a strict control Lyapunov function $V$ implies that
for some real number $r$ the vector $\sigma + r \tau$ is orthogonal to $DV$, so
$\tau\cdot DV\ne 0$ at all points where $\sigma \cdot DV\ne 0$, and we
can define for all $x\in \R^N\setminus \{0\}$
$$
h(x):=\left\{\begin{array}{ll} 0 & \text{ if } \sigma(x) \cdot DV(x) = 0 ,\\
- \frac{\sigma(x) \cdot DV(x)}{\tau(x) \cdot DV(x)} & \text{ if }\sigma(x) \cdot DV(x)\ne
0 . \end{array} \right.
$$
\begin{prop} Assume the system (\ref {affine2}) has a strict control
Lyapunov function $V\in C^2(\R^N\setminus\{0\})$ and the function $h$
is continuous in $\R^N\setminus\{0\}$.  Then there exists continuous functions
$k_i : \R^N\setminus\{0\} \to \R$, $i=1,...,P-1$, such that
$(k_1(x),...,k_{P-1}(x),h(x))$ is an almost surely asymptotically
stabilizing feedback for the system (\ref{affine2}).

Moreover, $k_i(x)\in [-1, 1]$ for $x$ in a neighborhood of 0 if $DV$
and $D^2V$ are bounded near 0, and either
$DV(x)\to 0$ or $g_i(x)\to 0$ for all $i$ as $x\to 0$.

\end{prop}
\begin{proof}
We recall from \cite{son1} that the function $\phi(a,0):=0$ for $a<0$,
$\phi(a,b):=(a+\sqrt{a^2+b^2})/b$ is real-analytic in the set
$S:=\{(a,b)\in\R^2\,:\, b>0 \text{ or } a<0\}$.  We set
\begin{multline*}
\gamma(x):= \\
f(x)\cdot DV(x) + trace\left[(\sigma(x)+h(x)\tau(x))(\sigma(x)+h(x)\tau(x))^T
\frac{D^2V(x)}{2}\right] +\frac{l(x)}{2} ,
\end{multline*}
$$
\beta(x) := \sum_{i=1}^{P-1} (g_i(x)\cdot DV(x))^2 .
$$
Since $V$ is a strict control Lyapunov function,
$$
\gamma(x) + \inf_{\alpha_i\in \R} \sum_{i=1}^{P-1}\alpha_i g_i(x) \cdot
DV(x)
\leq - \frac{l(x)}{2},
$$
so, for $x\ne 0$,
$$
\beta(x) = 0 \quad \Rightarrow \quad \gamma(x)\leq -l(x)/2 < 0 .
$$
Therefore $(\gamma(x),\beta(x))\in S$.
 Now we define, for $i=1,...,P-1$,
$$
k_i(x) := - \phi(\gamma(x),\beta(x)) g_i(x)\cdot DV(x), \quad x\ne 0,
$$
and $k(0)=0$.  Then $(k_1(x),...,k_{P-1}(x),h(x))$ is continuous in
$\R^N\setminus\{0\}$ and it satisfies
\begin{multline*} \left(f+\sum_{i=1}^{P-1}k_ig_i\right)\cdot DV +
trace\left[\left(\sigma+h\tau)(\sigma+h\tau\right)^T
\frac{D^2V}{2}\right] +\frac{l}{2} = \\
\gamma - \beta\phi(\gamma,\beta) = - \sqrt{\gamma^2+\beta^2} < 0 .
\end {multline*}
 Since $(\sigma+h\tau)\cdot DV = 0$ by definition of $h$, $V$ is a
strict Lyapunov function for the equation
$$
dX_t=\left(f(X_t)+\sum_{i=1}^{P-1} k_i (X_t) g_i(X_t)\right) dt+
\left(\sigma(X_t) + h(X_t) \tau(X_t) \right) dB_t .
$$
Therefore the origin is a.s.  asymptotically stable for this equation.

Finally we check the boundedness of $k$ in a neighborhood of 0.  This
is trivial for $\beta(x)=0$.  If $\beta(x)\ne 0$
$$
|k| \leq \frac{|\gamma + |\gamma| + \beta |}{\sqrt{\beta}}.
$$
Since either $DV\to0$ or $g_i\to 0$ for all $i$,
$\beta(x)\to 0$ as $x\to 0$.
We fix $\delta>0$ such that $\beta(x)\leq \delta$ implies
$\gamma(x)<0$, and then choose a neighborhood of the origin where
$\beta(x)\leq \delta$.  In this set $|k(x)| \leq \sqrt{\beta(x)}\to 0$.
\end{proof}
\begin{rem}\upshape
The proof above gives an explicit formula for the stabilizing feedback
in terms of the data and the Lyapunov function $V$ only, which reduces
to (\ref{formulafeed}) if $\tau\equiv 0$ and $P=2$.  From the formula
one sees that the feedback is $C^1$ in $\R^N\setminus\{0\}$ if $h, f,
g, \sigma, \tau$ and $l$ are such and $V\in C^3(\R^N\setminus\{0\})$.

Note also that the continuity assumption on $h$ is automatically
satisfied if $\tau\cdot DV$ is either always nonnull or identically $0$.

Finally, it is straightforward to extend the Proposition to the case
of $M$-dimensional noise with independent Brownian components
$B^1_t,..., B^M_t$ and diffusion term of the form
$
\sum_{i=P}^{P+M-1}\left(\sigma_i + \alpha^i_t \tau_i\right) dB^i_t ,
$
with $\sigma_i, \tau_i$ vector fields and $\alpha^i_t$ scalar controls.
\end{rem}
 %
\section{Some variants and extensions.}
In this section we collect several remarks on other applications of our methods.
We begin with the {\em almost sure exponential stabilizability}.
It means that there exists a positive rate $\lambda$ such that
for every initial data $x$ there exists an admissible control
$\oa_.\in {\cal A}_x$ whose corresponding trajectory $\ox_.$ satisfies
\[
V(\ox_t)\leq e^{-\lambda t} V(x) \quad \text{ a.s.}
\]
\begin{prop}
[a.s.  exponential stabilizability] Under the assumptions (\ref{condition1})
and (\ref{convex}), the null state is almost surely exponentially
stabilizable for (CSDE) if 
sufficient that the system admits there exists a control Lyapunov
function $V$ satisfying conditions (i), (ii), (iii) in the
Definition \ref{cliap} and for some $\lambda>0$

\noindent(iv)$^\prime$ for every $(p,Y)\in{\cal J}^{2,-}V(x)$
there exists $\overline{\alpha}\in A$ such that
\[
\sigma(x,\overline{\alpha})^T p=0 \; \,\mbox{ and }\, -p\cdot
f(x,\overline{\alpha})-trace\left[a(x,\overline{\alpha}) Y\right]
- \lambda V(x)\geq 0 .
\]
\end{prop}
\begin{proof}
We consider the system
\[
\left\{ \begin{array}{l} dX_t = f(X_t,\alpha_t)dt +
\sigma(X_t,\alpha_t)dB_t\\ dY_t = dt \end{array}\right.\] %
with initial data $X_0=x$ and $Y_0=0$, and the Lyapunov function
$W(x,y)=e^{\lambda y}V(x)$. By applying Theorem \ref{thmliap1} we obtain
the existence of a control $\oa_.$ such that the corresponding trajectory
almost surely satisfies $V(\ox_t)\leq V(x)e^{-\lambda t}$, which is
the desired inequality.
\end{proof}
%
%
%

Next we extend the results of Section 2 to the stabilizability of
a general closed set $M\subseteq\R^N$.  We denote with $d(x,M)$ the distance
between a point $x\in \R^N$ and $M$.
 \begin{defi}[a.s.  stabilizability  at $M$]\label{def:mstab}
The system $(CSDE)$ is \emph{almost surely
(open loop) stabilizable} at $M$ if 
there exists $\gamma\in\mathcal{K}$ such
that, for every $x$ in a neighborhood of $M$, there is an admissible control function $\overline{\alpha}_{\cdot}\in{\cal A}_x$
whose trajectory
 $\overline{X}_{\cdot}
$
verifies
$$
d(\overline{X}_t,M)\leq \gamma(d(x,M)) \quad
\forall \, t\geq 0 \quad
\text{almost surely.}
$$

If, in addition,
$$
\lim_{t\rightarrow +\infty}d(\ox_t,M)=0 \quad a.s.
$$
the system is  \emph{almost surely (open loop) locally asymptotically
stabilizable} at $M$.

If these properties hold for all $x\in\R^N$ the system is  \emph{a. s. (open loop)
globally asymptotically stabilizable} at $M$.
\end{defi}

\begin{rem}\upshape
If $M$ is a.s.  stabilizable, then it is viable for $(CSDE)$.  In
fact, the definition
%
gives for $x\in M$ and $\eps>0$
an admissible control such that
almost surely
$d(X_t,M)\leq \eps$ for all $t\geq 0$.
Then for such control and any $\lambda>0$ almost surely
$\mathbf{E}_x\int_0^{+\infty} d(X_t,M)e^{-\lambda t}dt \leq
\frac{\eps}{\lambda}$
and so
\[
\inf_{\alpha_.\in {\cal A}_x}  \mathbf{E}_x\int_0^{+\infty}
d(X_t,M)
e^{-\lambda t}dt=0 .\]
The
convexity assumption (\ref{convex}) and an existence theorem for
optimal controls \cite{hl} imply that the inf is attained,
and the minimizing control produces a trajectory staying in $M$
for all $t\geq 0$.
\end{rem}

\begin{defi}[control Lyapunov functions at $M$]
\label{cmliap} Let ${\cal O}$ be an open neighborhood of the
closed set $M$. A function $V : {\cal O} \rightarrow[0,+\infty)$ is
a control Lyapunov function at $M$ for $(CSDE)$ if

\noindent (i) $V$  is lower semicontinuous;

\noindent(ii) 
there exists $\gamma_1\in\mathcal{K}_\infty$ such that
$
V(x)\leq\gamma_1 (d(x,M))
$
for all $x\in{\cal O}$;

\noindent(iii) 
there exists $\gamma_2\in\mathcal{K}_{\infty}$ such that
$\gamma_2(d(x,M))\leq V(x)$ for all $x\in{\cal O}$;


\noindent (iv) for all $x\in{\cal O}\setminus M$ and
$(p,Y)\in{\cal J}^{2,-}V(x)$ there exists $\overline{\alpha}\in A$
such that condition (\ref{ceq}) holds.

The function $V$ is a \emph{strict control Lyapunov function} at $M$
if it satisfies conditions
(i)-(iii)
and

\noindent(iv)$^{\prime}$ for some Lipschitz continuous $l : {\cal
O}\to \R$, $l(x)>0$ for all $x\notin M$, and $(p,Y)\in{\cal
J}^{2,-}V(x)$ there exists $\overline{\alpha}\in A$ such that
condition (\ref{cseq}) holds.

\end{defi}
%
%
%

Now we can state the analogues of the First and Second Lyapunov Theorems
for the a.s. stabilizability at $M$. Their proofs are easily obtained
from the arguments of the Theorems \ref{thmliap1} and \ref{thmliap2} by using $d(x,M)$
instead of $|x|$ and noting that
 conditions $(ii)$ and $(iii)$ in the Definition \ref{cmliap}
say that the sublevel sets of the Lyapunov function
form a basis of neighborhoods of $M$.

\begin{teo}\label{Mthmliap1}
Assume (\ref{condition1}), (\ref{convex}), and the existence of a control
Lyapunov function $V$ at $M$.  Then

\noindent (i) the system $(CSDE)$ is almost surely stabilizable at $M$;

\noindent (ii) if, in addition, the domain ${\cal O}$ of $V$ is all $\R^N$,
for all $x\notin M$ there
exists $\oa_.\in {\cal A}_x$ such that the corresponding trajectory
$\ox_.$ satisfies
$$
d(\ox_t,M) \leq \gamma_1^{-1}(\gamma_2(dist(x,M)))\quad \forall\,
t\geq 0\quad{a.s.}
$$
with $\gamma_1, \gamma_2\in {\cal K}_{\infty}$ from Definition \ref{cmliap};
in particular, if $M$ is bounded, the system is also a.s.  Lagrange stabilizable.
\end{teo}
\begin{teo}
    \label{Mthmliap2}
Assume (\ref{condition1}), (\ref{convex}), and the existence of a
strict control Lyapunov function $V$ at $M$.  Then

\noindent (i) the system $(CSDE)$ is a.s. locally asymptotically stabilizable
at $M$;

\noindent (ii) if, in addition, the domain ${\cal O}$ of $V$ is all $\R^N$,
the system is a.s. globally asymptotically stabilizable at $M$.
\end{teo}

\begin{rem}[Stochastic target problems and absorbing sets]\upshape
%
A stochastic target problem
consists of steering the state of the system $(CSDE)$ in finite
time into a given closed set $\mathcal{T}$ (the target) by an
appropriate choice of the control.  One of the objects of interest is
the set of initial positions from which this goal can be achieved
almost surely in a given time $t$.  We define these reachability sets
for $t>0$ as \[\mathcal{R}(t)=\{x\in\R^N\ |\ \exists\alpha_{.}\in
{\cal A}_x\,: \; X_{t}\in\mathcal{T} \mbox{ a.  s.} \}\] We consider a
target  $\mathcal{T}$ 
containing $0$ and invariant for the stochastic
system and we assume there exists a global strict control Lyapunov
function $V$ as defined in (\ref{scliap}) such that
\[\inf_{\R^N\setminus \mathcal{T}}l(x)=L>0 .\]
We are going to show that each reachability set $\mathcal{R}(t)$ lies between
two sublevel sets of the Lyapunov function $V$.
The arguments in the proof of Theorem \ref{thmliap2},
show that for every initial point
$x\not\in \mathcal{T}$ there exists a control
$\oa{.}\in {\cal A}_x$ such that the first entry time
$\overline{\tau}_x$ of the corresponding trajectory in the target
 is almost surely bounded by
\begin{equation}\label{target}
\overline{\tau}_x\leq \left(
V(x)-\inf_{\partial \mathcal{T}}V(y)\right)/L .
\end{equation}
In particular, since the target $\mathcal{T}$ is invariant, it is
reached almost surely in a finite time, and such time is also uniformly
bounded, $\mathcal{T}$ is an {\em absorbing set} for the system
according to the terminology in ~\cite{ad4}.  %
Next, from the assumptions and inequality (\ref{target}) we get
\[
\{x\in\R^N\ |\ V(x)\leq Lt+\inf_{\partial
\mathcal{T}}V(y)\} \subseteq \mathcal{R}(t) .
\]
Using Chebyshev
inequality and estimates of the same kind as in the proof of Theorem
\ref{thmliap2} we can find also for every $t>0$ a positive number
$k(t)$ depending continuously on $t$ such that
\[
\mathcal{R}(t)\subseteq \{x\in\R^N\ |\ V(x)\leq k(t)\}.
\]
Let us mention that Soner and Touzi ~\cite{st1} developed  recently a
PDE approach to stochastic target problems, see also ~\cite{st2} and
the references therein for some interesting applications to geometric
PDEs and front propagation problems.  %
\end{rem}

\section{Examples}

We begin with an example of an uncontrolled system that does not have a
continuous Lyapunov function but has a l.s.c.  Lyapunov function and
therefore is a.s. stable.  It shows that allowing $V$ to be merely l.s.c.
in Theorem \ref{thmliap1} really increases the range of the
applications.  Our example is a variant of a deterministic one by
Krasovskii \cite{kra}, namely,
 \[
 \left\{
\begin{array}{l}
\dot X_t=Y_t\\ \dot Y_t=-X_t+
Y_t(X_t^2+Y_t^2)^3\sin^2\left(\frac{\pi}{X_t^2+Y_t^2}\right) ,
\end{array}\right.
\]
see \cite{bacr} for a discussion of this and other deterministic examples.

 \begin{ex}\label{kras}\upshape 
%
We transform the previous system in polar coordinates and perturb it with a
white noise tangential to the circles
$C_n:=\{(x,y)\,:\, |(x,y)|=\frac{1}{\sqrt{n}}\}$ and
nondegenerate  between two consecutive circles:
\[
 \left\{
\begin{array}{lll}
d\rho_t & = & \left[ \rho_t^7
\sin^2(\theta_t)\sin^2(\frac{\pi}{\rho_t^2})\right] dt +
\left[\sigma(\rho_t,\theta_t)\sin^2(\frac{\pi}{\rho_t^2})\right]dB_t\\
d\theta_t & = &
\left[-1+\rho_t^6\sin(\theta_t)\cos(\theta_t)\sin^2(\frac{\pi}{\rho_t^2})\right]dt,
\end{array}\right.
\]
where $B_t$ is a 1-dimensional Brownian motion and $\sigma$
satisfies the hypotheses for the existence and uniqueness of the
solution of the stochastic differential equation.  As in the
undisturbed case, the circles $C_n$ are a.s.  invariant and any point
in $C_n$ is eventually reached a.s.  by any trajectory starting in
$C_n$.  Then any Lyapunov function $V$ is constant on $C_n$ because
$V(\rho_t,\theta_t)\leq V(\rho_0,\theta_0)$ a.s., and
$c_n:=V_{|C_n}\ne c_{n-1}:=V_{|C_{n-1}}$
at least on a subsequence.  By property $(iv)$ in the
definition \ref{cliap} of Lyapunov function,
for every $(\rho, \theta)$ in the interior of $C_{n-1}\setminus C_n$
and every $(p,X)\in {\cal J}^{2,-}V(\rho,\theta)$, we get
$\left(\sigma(\rho,\theta)\sin^2(\frac{\pi}{\rho^2}),0\right)\cdot p=0$.
Since the diffusion is nondegenerate in the $\rho$
direction in the interior of
$C_{n-1}\setminus C_n$, from the previous equality we deduce
 that, for such  $(\rho, \theta)$, every element in ${
\cal J}^{2,-}V(\rho,\theta)$  is of the form $\left((0,p_2), X\right)$.
This  implies that the function $V$ is  constant in the $\rho$
direction in the  interior of
$C_{n-1}\setminus C_n$ and then it cannot be continuous.

Now we check that the Lyapunov function of the undisturbed
system in the unit ball does the job also for our perturbed
stochastic system. We take
$$
V(\rho,\theta):=\frac{1}{\sqrt{n}}\quad\text{ for }
\frac{1}{\sqrt{n}}<\rho\leq\frac{1}{\sqrt{n-1}}, \quad\forall \theta.
$$
This is a positive definite function, lower semicontinuous and
continuous at $0$.  We calculate its second order subjets and plug them
into (\ref{ceq}).  If $\rho\ne\frac{1}{\sqrt{n}}$
for all $n$, $(p, \mathbf{X})\in J^{2,-}V(\rho,\theta)$ if and only if
$p=0$ and $\mathbf{X}\leq 0$, so the condition (\ref{ceq})
is trivially satisfied.
%

On the other hand,  $(p, \mathbf{X})\in
J^{2,-}V(\frac{1}{\sqrt{n}},\theta)$ if and only if
$$
p = \left(\begin{array}{c}s\\ 0 \end{array}\right), \; s\geq 0 \quad
\text{ and } \quad \mathbf{X} = \left(\begin{array}{cc}a & b\\ b &
c\end{array}\right), \; c\leq 0.
$$
%
At the points with $\rho=\frac{1}{\sqrt{n}}$ the drift $f$ of the system
is $(0, -1)$ and the dispersion
vector $\sigma$ is $(0,0)$.  Then
$$
f\cdot p + \frac{1}{2} trace\left[\sigma\sigma^T \mathbf{X}\right] =
0 , \qquad \sigma \cdot p = 0
$$
and the condition (\ref{ceq}) is satisfied.
%
Therefore Theorem \ref{thmliap1} applies and the system is a.s.  Lyapunov
stable at the origin.
  %
\end{ex}

\medskip
 The next two examples are about {\em stochastic perturbations of
stabilizable systems}.  We consider a deterministic controlled system
in $\R^N$
\begin{equation}
\label{detsyst}
\dot X_t=f(X_t,\alpha_t) 
\end{equation}
globally asymptotically (open loop) stabilizable at the origin, i.e.,
asymptotically controllable in the  terminology of
deterministic systems \cite{son, sosu}.  By the converse Lyapunov
theorem of Sontag \cite{son0, sosu}, there exists a strict continuous control
Lyapunov function for the system, i.e., for some positive definite
continuous function $L$, a proper function $V$ satisfying in
$\R^N\setminus\{0\}$
 \begin{equation}
 \label{ldet}
  \max_{\alpha\in A}\left\{ -f(x,\alpha)\cdot DV\right\} - L(x)\geq 0
  \end {equation}
 in the viscosity sense.  (This is perhaps not explicitly stated in the
 literature; the original result of Sontag \cite{son0} interprets this
 inequality in the sense of Dini derivatives of $V$ along relaxed
 trajectories, the paper of Sontag and Sussmann ~\cite{sosu} in the
 sense of directional Dini subderivatives, and both these senses are
 known to be equivalent to the viscosity one, see, e.g., ~\cite{subb, bcd}).
%
%
In the following examples we perturb in two different ways
(\ref{detsyst}) and give condition under which $V$ remains a control
Lyapunov function for the a.s.  stabilizability of the new stochastic system.
\begin{ex}\upshape
Consider the controlled diffusion process
\begin{equation}
\label{pertsyst1}
dX_t=f(X_t,\alpha)dt+\sigma(X_t)dB_t
\end{equation}
where $B_t$ is a $M$-dimensional Brownian motion and $\sigma$ a
 Lipschitzean $N\times M$ matrix.  Then

 \noindent{\em $V$ is a Lyapunov function for (\ref{pertsyst1}) if, for
 some open set ${\cal O}\ni 0$ and some continuous $l : {\cal O}\to
 [0,+\infty)$, $V$ satisfies in viscosity sense in ${\cal O}\setminus
 \{0\}$
 \begin{equation}
 \label{es2}
 -trace\left[\frac{1}{2}\sigma \sigma^T D^2V\right] + L - l \geq 0, \quad
 \sigma_i \cdot DV=0 \quad \forall i,
 \end{equation}
  and it is a strict Lyapunov function if $l$
 is positive definite. }

 In fact, this inequality and (\ref{ldet}) give,
 for any $(p,\mathbf{X})\in J^{2,-}V(x)$,
 $$
 \max_{\alpha\in A}\left\{ -f(x,\alpha)\cdot
 p\right\}-trace\left[\frac{1}{2}\sigma \sigma^T \mathbf{X}\right] - l \geq 0,
 $$
so $V$ satisfies the inequality  in condition (\ref{cseq}),
whereas the equality in condition (\ref{cseq}) reduces
to $\sigma_i\cdot p=0$.
 %

In the classical special case of $V(x)=|x|^2$ and $M=1$, the
sufficient condition (\ref{es2}) for $V$ to be a Lyapunov function of
(\ref{pertsyst1}) reads
$$
l(x):=L(x) - |\sigma(x)|^2\geq 0, \quad \sigma(x)\cdot x=0.
$$
For a noise of dimension $M=N$ an example of $\sigma$ satisfying the
orthogonality condition in (\ref{es2}) is $\sigma(x)=k
 \left(\mathbf{I}-\frac{DV(x)\otimes DV(x)}{|DV(x)|^2} \right)$ for any
 constant $k$.
\end{ex}
\begin{ex}\upshape
Here we consider the perturbation of the deterministic system
(\ref{detsyst}) by a function $g$ of a $K$-dimensional diffusion
process $Y_t$:
\begin{equation}
\label{pertsyst2}
\left\{\begin{array}{l} \dot X_t =
f(X_t,\alpha_t) + g(X_t,Y_t)\\ dY_t=b(Y_t, X_t, \alpha_t)dt
+\tau(Y_t, X_t, \alpha_t)dB_t
\end{array}
\right.
\end {equation}
 where the function $g:\R^n\times\R^K\to \R^n$ is Lipschitz
continuous
%
%
with $g(0,y)=0$ for all $y$, $B_t$ is a 1-dimensional Brownian motion,
and $b, \tau$ are vector fields in $\R^K$ with the usual assumptions.  We
are still assuming that (\ref{detsyst}) has a strict control Lyapunov
function $V$, i.e., (\ref{ldet}) holds with $L$ positive definite.  We
are interested in the stabilizability of the perturbed system at the
set $M:=\{(x,y)\in \R^n\times\R^K \,:\; x=0\}$, which corresponds to
the origin of the unperturbed system (\ref{detsyst}), see Definition
\ref{def:mstab}.  Note that the assumption on $g$ implies the
viability of $M$ for (\ref{pertsyst2}).  We claim that

 \noindent{\em the function $V$, defined by $V(x,y):=V(x)$ for all $y$, is a
 Lyapunov function at $M$ for (\ref{pertsyst2}) (see Definition
 \ref{cmliap}) if, for some open set ${\cal O}\ni 0$ and some
 continuous $l : {\cal O}\times\R^K \to [0,+\infty)$,
$V$ satisfies in
 viscosity sense in ${\cal O}\setminus\{0\}$
\begin{equation}
 \label{es3}
\inf_{y\in\R^K}\{-g(x,y)\cdot DV(x) - l(x,y)\} + L(x) \geq 0,
 \end{equation}
  and $V$ is a strict Lyapunov function if $l(x,y)>0$
 for all $x\ne 0$ and all $y$.  }

In fact, since $d((x,y),M)=|x|$, $V$ satisfies the conditions
$(i)-(iii)$ of the Definition \ref{cmliap}.
By (\ref{ldet}) and (\ref{es3})
$V$ is also a viscosity supersolution in ${\cal O}\times \R^K\setminus M$ of
\[
\sup_{a\in A}\left\{- f(x,a)\cdot DV(x)\right\}-g(x,y)\cdot DV(x) - l(x,y)
\geq  0,
\]
which is the inequality in (\ref{cseq}) in this case, because $V$ is
constant in $y$.  Finally, for the same reason, the condition in
(\ref{cseq}) of orthogonality of the diffusion vector to the level
sets of $V$ is trivially satisfied.

%

The inequality (\ref{es3}) is a smallness condition of the component
of $g$ in the direction of $DV$ with respect to $L$ in the set ${\cal
O}$, uniformly in $y$.  For $l\equiv 0$ and $V$ smooth in ${\cal
O}\setminus\{0\}$ it becomes
\begin{equation}
 \label{es3bis}
\sup_{y\in\R^K}g(x,y)\cdot DV(x)\leq L(x), \quad {in }\; {\cal O}\setminus\{0\},
\end {equation}
which is satisfied, in particular, if
$$
\sup_{y\in\R^K}|g(x,y)|\leq L(x)/ Lip V,
$$
where $Lip V$ denotes the Lipschitz constant of $V$ in ${\cal O}$.
We recall that, under our assumption that the deterministic system
(\ref{detsyst}) be asymptotically controllable, although $V$ may not
be smooth, it can be chosen semiconcave in $R^N\setminus\{0\}$ and
therefore locally Lipschitz ~\cite{r}.  If we make this choice, it is
enough that the inequality (\ref{es3bis}) holds for all points %
$x\in{\cal O}$
where $V$ is differentiable, and the last inequality is guaranteed for
all perturbations $g$ with small sup-norm with respect to $y$.

\end{ex}

\medskip
In the next two examples we give conditions on a radial function to be
a Lyapunov function for a.s.  stability.
\begin{ex}\label{exrad}\upshape
We consider as a candidate Lyapunov function for the general
controlled system $(CSDE)$ the function $V(x)=v(|x|)$, for some smooth
$v:[0,+\infty)\to [0,+\infty)$ with $v'(r)>0$ for $r>0$.  Since
$DV(x)=xv'(|x|)/|x|$, in view of the orthogonality condition in
(\ref{ceq}) we restrict ourselves to controls $\alpha \in A$ such
that
\begin{equation}
\label{ortrad}
\sigma_i(x,\alpha)\cdot x = 0, \quad \forall\, i= 1,\ldots, M.
\end{equation}
We compute
$$
trace\left[a(x,\alpha) D^2 V(x)\right] = \frac{v'(|x|)}{|x|}
trace\  a(x,\alpha) + \left(v''(|x|) - \frac{v'(|x|)}{|x|}\right)
\frac{|\sigma(x,\alpha)^T x|^2}{|x|^2}
$$
and use (\ref{ortrad}) to obtain that  $V$ is a Lyapunov function
if and only if, in a neighborhood ${\cal O}$ of 0,
$$
l(x):=\max_{\alpha\in A,\, \sigma(x,\alpha)^T x = 0}
\left[-f(x,\alpha)\cdot x - trace\, a(x,\alpha)
\right]\frac{v'(|x|)}{|x|} \geq 0,
$$
i.e.,
\begin{equation}
\label{condrad}
\min_{\alpha\in A,\, \sigma(x,\alpha)^T x = 0}\left[f(x,\alpha)\cdot x + trace\, a(x,\alpha)
\right] \leq 0.
\end{equation}
This condition is independent of the choice of $v$.  Moreover, if $l>0$ and
Lipschitz in  ${\cal O}\setminus\{0\}$  and $l\to 0$ as $x\to 0$,
then $V$ is a strict Lyapunov function. Note that, although the
radial component of the diffusion must be null by (\ref{ortrad}), its
rotational component still plays a destabilizing role. In fact,
$trace\, a(x,\alpha)\geq 0$ and whenever it is nonnull it must be
compensated by a negative radial component of $f$.

In particular, a single-input affine system with uncontrolled diffusion and
1-dimensional noise $B_t$
\begin{equation*}
 dX_t=\left(f(X_t)+\alpha_t g(X_t)\right)dt+
\sigma(X_t) dB_t , \quad \alpha_t\in[-1, 1],
\end{equation*}
has a radial Lyapunov function in ${\cal O}$ if and only if
$$
\sigma(x)\cdot x=0 \quad\text{and}\quad |g(x)\cdot x|\geq f(x)\cdot
x + \frac{|\sigma(x)|^2}{2}\quad \text{in }{\cal O},
$$
and $V(x)=|x|^2/2$ is a strict Lyapunov function in ${\cal O}$ if and only if
$$
l(x):= |g(x)\cdot x| - f(x)\cdot
x - \frac{|\sigma(x)|^2}{2}>0 \quad \text{in }{\cal O}\setminus\{0\}.
$$
Moreover, $k(x):=-sign(g(x)\cdot x)$ is a stabilizing feedback if
$g(x)\cdot x$ does not change sign; if it does $k$ is discontinuous,
then a
continuous stabilizing feedback in a neighborhood of $0$ is given by
the formula (\ref{formulafeed}) in Section 4.
\end{ex}
\begin{ex}\upshape
Here we study a system in $\R^2$ written in polar coordinates $(\rho,
\theta)$ and look for radial Lyapunov functions, i.e., of the form
$V(\rho,\theta)=v(\rho)$.  Consider the stochastic controlled system:
\[
(CSDE)\left\{\begin{array}{l}
d\rho_t=f(\rho_t,\theta_t,\alpha_t)dt+
\sigma(\rho_t,\theta_t,\alpha_t)dB_t\\
d\theta_t=g(\rho_t,\theta_t,\alpha)dt+
\tau(\rho_t,\theta_t,\alpha_t)dB_t
\end{array}
\right.
\]
where
all functions $f, \sigma, g, \tau$ are $2\pi$-periodic and $B_t$ is (for
simplicity) a $1$-dimensional Brownian motion.
%
The conditions for a function $V=v(\rho)$ to be a Lyapunov function of
this system at the set $M:=\{(0,\theta)\,:\, \theta\in\R\}$ are the
following.  The orthogonality condition in (\ref{cseq}) requires that
for every $(\rho,\theta)$ there exists a subset $A(\rho,
\theta)\neq\emptyset$ of the control set $A$ such that
\[
\sigma(\rho,\theta,\alpha)=0\ \ \ \ \ \forall \alpha\in A(\rho, \theta).
\]
Then the condition (\ref{cseq}) is satisfied if
$v$ is a viscosity supersolution of the ordinary differential inequality
 \[
 \sup_{\alpha\in A(\rho, \theta)}\left\{ -f(\rho,\theta,\alpha)\cdot v'(\rho)
\right\}\geq 0
\]
for $\rho>0$ and for each fixed $\theta\in[0,2\pi]$.  Of course the
same result can be obtained from the previous example with some
calculations based on the Ito chain rule.
\end{ex}
\medskip
The last two examples  are about the stabilization to sets $M$ different
from the origin,
namely, the complement of a ball and a periodic orbit.

\begin{ex}\upshape \label{exexterior}
%
%
We consider
the general system $(CSDE)$ and the set
$$
M := \{x\  |\ |x|\geq R\} = \R^N\setminus B_R.
$$
We assume $M$ is viable for the system. We
take the radial function $V$ 
\[
V(x):=\left\{\begin{array}{lll}R^2-|x|^2& & |x|<R \\ 0& & |x|\geq
R\end{array}
\right.
\]
and use the calculations of Example \ref{exrad} to see that $V$ is a
Lyapunov function at $M$ if and only if
%
for every $x$ with $|x|<R$
there exists $\oa\in A$ such that
\[
\sigma_i(x,\oa)\cdot x=0 \quad \forall i
\quad\text{and}\quad f(x,\oa)\cdot x + trace\  a(x,\oa)\geq 0 .
\]
On the contrary of Example \ref{exrad}, here the rotational component
of the diffusion has a stabilizing effect. In fact, the drift $f(x,a)$
is allowed also to point away from $M$ if its negative radial component
is compensated by the positive term $trace\  a(x,\oa)$.

If $K\subset B_R$ is a compact set and
$$
l(x):=\max_{\alpha\in A,\, \sigma(x,\alpha)^T x = 0}\left[f(x,\alpha)\cdot x + trace\, a(x,\alpha)
\right] > 0 \quad \text{in }B_R\setminus K ,
$$
then $M$ is locally asymptotically stable by Theorem \ref{Mthmliap2} and for all
 initial points $x\notin K$ there is a control whose trajectories tend a.s. to $M$ as
$t\to+\infty$. In this case we can say that $K$ can be made almost surely
repulsive by a suitable choice of the controls. In particular, we have
a criterion of instability of an equilibrium point.

Note also that if $l>0$ on $\partial M = \partial B_R$ then for  some
control the
trajectories starting in a suitable neighborhood of $\partial M$
reach $M$ in finite time a.s., as we observed in the last remark of
Section 5. In particular, if $l>0$ in
$\overline{B_R}$, then for every $x\in B_R$  there
exists a control $\oa_.$ such that
the exit time of the corresponding trajectory $\ox_.$ from $B_R$
is almost surely bounded by $(R^2-|x|^2)/\min_{B_R} l$.

\end{ex}

\begin{ex}\upshape
Consider $(CSDE)$ in $\R^2$ and assume the circle $\gamma:=\{ x\,:\,
|x|=R \}$ is a viable set.  By the results of \cite{bj} this occurs if
for all $x\in\gamma$ there exists $\oa\in A$ such that
$$
\sigma(x,\oa)\cdot x=0\quad\text{ and }\quad f(x,\oa)\cdot x + trace\,
a(x,\oa)=0 .
$$

Then $\gamma$ is
locally asymptotically stabilizable if,
in a neighborhood $\{ x\,:\, R-\varepsilon \leq |x|\leq R +\varepsilon \}$,
\[
\max_{\alpha\in A,\, \sigma(x,\alpha)^T x = 0}\left[ f(x,\alpha)\cdot x +
trace\  a(x,\alpha)\right] > 0 \quad \text{ if }\, |x|<R ,
\]
\[
\min_{\alpha\in A,\, \sigma(x,\alpha)^T x = 0}\left[ f(x,\alpha)\cdot x +
trace\  a(x,\alpha)\right] < 0 \quad \text{ if }\, |x|>R .
\]
This follows immediately from the arguments of the Examples
\ref{exrad} and \ref{exexterior}.
\end{ex}


\end{document}